\documentclass{siamart220329}

\usepackage{tim,epsfig,fullpage,algorithm,algorithmic}

\usepackage{cite}
\usepackage{ma580}
\usepackage{graphicx}
\usepackage{hyperref}
\usepackage{xcolor,mdframed}
\usepackage{url}
\usepackage{doi}
\def\calf{{\cal F}}

\def\pkg{{\bf MultiPrecisionArrays.jl} }
\def\version{{\ v0.1.4 \ }}

\usepackage{listings}
\usepackage[nocolors]{jlcode}
\lstdefinestyle{code}{
  language=Julia,
  xleftmargin=1.5cm,
  showstringspaces=false,
}
\lstnewenvironment{code}{
\lstset{
  language=Julia,
  xleftmargin=1.5cm,
  showstringspaces=false,
}}{}

\title{Using MultiPrecisonArrays.jl: \\ Iterative Refinement in Julia}
\author{
C. T. Kelley%
\thanks{North Carolina State University,
Department of Mathematics,
Box 8205, Raleigh, NC 27695-8205, USA
(ctk@ncsu.edu).
This work was partially supported by
Department of Energy grant DE-NA003967.
}
}

\begin{document}
\maketitle

\begin{abstract}
MultiPrecisionArrays.jl is a Julia package.
This package provides data structures and solvers for several 
variants of iterative refinement. It will become much more useful 
when half precision (aka Float16) is fully supported in LAPACK/BLAS. 
For now, its best general-purpose application is classical iterative refinement with double precision equations and single precision factorizations.

It is useful as it stands for people doing research in iterative refinement.
We provide a half precision LU factorization that, while far from 
optimal, is much better than the default in Julia.

This document is for \version of the package. We used both the
AppleAccelerate \cite{accelerate} framework and openBLAS for the examples 
and the timings may be different from the older versions. 
All the computations in this paper were done on an Apple Mac Mini with
an M2 Pro processor and 8 performance cores. The examples were done with
several different versions of \pkg, Julia, and the OS.
Hence the examples are indicative of the performance
of IR, but the reader's results may differ in timings and iteration counts.
The qualitative results are the same.
\end{abstract}



\section{Introduction}
\label{sec:intro}

The Julia \cite{Juliasirev}
package {\bf MultiPrecisionArrays.jl} \cite{ctk:mparrays,ctk:joss2}
provides data
structures and algorithms for several variations of iterative 
refinement (IR). In this introductory section we look at the 
classic version of iterative refinement and discuss its implementation 
and convergence properties.

We focus on two views of IR. One is as a
perfect example of a storage/time tradeoff.
To solve a linear system $\ma \vx = \vb$ in $R^N$ with IR,
one incurs the storage penalty of making a low
precision copy of $\ma$ and reaps the benefit of only having to
factor the low precision copy. 

The other viewpoint is as a way to address 
very ill-conditioned problems, which was the purpose of the 
original paper on the topic \cite{Wilkinson48}. 
Here one stores and factors $\ma$
in one precision but evaluates the residual 
in a higher precision. 

In most of this paper we consider IR using two precisions. We use the
naming convention of 
\cite{amestoy:2024,CarsonHigham1,CarsonHigham,demmelir} for these
precisions and the others we will consider later. We start with the
case where we store the data $\ma$ and $\vb$ in the higher or {\em working}
precision and compute the factorization in the lower or {\em factorization}
precision. Following standard
Julia type notation, we will let $TW$ and $TF$ be the working and factorization
precision types. So, for example
\begin{code}
x = zeros(TW,N)
\end{code}
is a high precision vector of length $N$.

In our first case, for example, $TW$ will be double
and $TF$ will be single. In the second case, $TW$ and $TF$
will both be single precision and we will compute the residual
in precision $TR$, which will be double.
We will make precision and inter precision
transfers explicit in our algorithmic descriptions. 

The first three sections of this paper use \pkg to generate 
tables which compare the algorithmic options, but do not talk about
using Julia.
Algorithm~\ref{alg:ir} is
the textbook version \cite{higham} version of the algorithm for the
$LU$ factorization.

\begin{algorithm}
\label{alg:ir}
{$\mbox{\bf IR}(\ma, \vb)$}
\begin{algorithmic}
\STATE $\vx = 0$
\STATE $\vr = \vb$
\STATE Factor $\ma = \ml \mU$ in precision $TF$
\WHILE{$\| \vr \|$ too large}
\STATE $\vd = \mU^{-1} \ml^{-1} \vr$
\STATE $\vx \leftarrow \vx + \vd$
\STATE Evaluate $\vr = \vb - \ma \vx$ in precision $TR$
\ENDWHILE
\end{algorithmic}
\end{algorithm}

One must be clear on the meanings of ``factor in precision $TF$'' and
$\vd = \mU^{-1} \ml^{-1} \vr$ to implement the algorithm. As we indicated
above, if $TW$ is a higher precision than $TF$, 
the only way to factor $\ma$ in low precision is to make a copy
and factor that copy. We must introduce some notation for that. We 
let $\calf_p$ be the set of floating point numbers in precision $p$,
$u_p$ the unit roundoff in that precision, and
$fl_p$ the rounding operator. Similarly we let
$\calf^N_p$, $\calf^{N \times N}_p$ denote the vectors and matrices
in precision $p$. We let $I_p^q$ denote the copying operator from
precision $p$ to precision $q$. When we are not explicitly specifying
the precisions, we will use $W$ and $F$ as sub and superscripts for
the working and factorization precisions. 
When we are discussing specific use cases
out sub and superscripts will the $d$, $s$, and $h$ for double, single,
and half precision.

So, factoring a high-precision matrix $\ma \in \calf^{N \times N}_H$ 
in low precision $TF$ means copy $\ma$ into low precision and obtain
\[
\ma_F = I_W^F (\ma) 
\]
and then factor $\ma_F = \ml \mU$. 

The current version of \pkg \version requires that the type of $\ma$ be \\
{\tt AbstractArray\{TW,2\}} where $TW$ is single or double. 
We will make this
more general in a later version, but we will always require that the
Julia function $lu!$ accept the type of $\ma$. In particular, this means
that \pkg will not accept sparse arrays. I'd like to fix this, but 
have no idea how to do it.

\section{Integral Equations Example}
\label{sec:integeq}

The submodule {\bf MultiPrecisionArrays.Examples} has an example which
we will use repeatedly. The function {\bf Gmat(N)} returns the $N$ point
trapezoid rule discretization of the Greens operator
for $-d^2/dx^2$ on $[0,1]$

\[
G u(x) = \int_0^1 g(x,y) u(y) \, dy
\]

where

\[
g(x,y) =
    \left\{\begin{array}{c}
        y (1-x) ; \ x > y\\
        x (1-y) ; \ x \le y
    \end{array}\right.
\]

The eigenvalues of $G$ are $1/(n^2 \pi^2)$ for $n = 1, 2, \dots$.

The code for this is in the {\bf /src/Examples} directory.
The file is {\bf Gmat.jl}.

In the examples we will use {\bf Gmat} to build a matrix
$\ma = I - \alpha \mg$. In the examples
we will use $\alpha=1.0$, a very well conditioned case, 
$\alpha=800.0$, a somewhat ill-conditioned case, 
and $\alpha=799.0$, a very ill-conditioned case.
The latter case is very near singularity.

\subsection{Two Classic Examples}
\label{sec:classic}
While \pkg was designed for research,
it is useful in applications. We present two examples. 
The first is one where $TW = TR$ is double and $TF$ is single. 
This example illustrates the essential ideas in IR and shows
how to use \pkg in the most direct way. We follow this example with
an illustration on how one harvests iteration statistics.

The second example is more subtle. Here $TR$ is double, $TF=TW$
is single, and the matrix is highly ill-conditioned. Here the termination
criteria is different and we look at more iteration statistics.

\subsubsection{Double-Single Precision}
\label{subsubsec:doublesingle}
In this case the working precision is double 
and the factorization precision is single.
This case avoids the (very interesting) problems with half precision. 

Here is a Julia code that implements IR in this case. We will use
this as motivation for the data structures in \pkg.

\begin{code}
""" 
IR(A,b)
Simple minded iterative refinement
Solve Ax=b
"""
function IR(A, b)
    x = zeros(length(b))
    r = copy(b)
    tol = 10.0 * eps(Float64)
    #
    # Allocate a single precision copy of A and factor in place
    #
    A32 = Float32.(A)
    AF = lu!(A32)
    #
    # Give IR at most ten iterations, which it should not need
    # in this case
    #
    itcount = 0
    rnorm=norm(r)
    rnormold = 2.0*rnorm
    while (rnorm > tol * norm(b)) && (rnorm < .9 * rnormold)
        #
        # Store r and d = AF\r in the same place.
        #
        ldiv!(AF, r)
        x .+= r
        r .= b - A * x
        rnorm=norm(r)
        itcount += 1
    end
    return x
end
\end{code}

\subsection{Running MultiprecisionArrays: I}
\label{sec:running1}

The function {\bf IR}  allocates memory for the low
precision matrix and the residual with each call. \pkg addresses
that with the {\bf MPArray} data structure which
allocates for the low precision copy of $\ma$ and the residual $\vr$.

The most simple way to use this package is to combine the construction
of the MPArray with the factorization of the low precision copy of $\ma$.
One does this with the {\bf mplu} command.

As an example we will solve the integral equation 
with both double precision $LU$ and an MPArray
and compare execution time and the quality of the results.
We will use the function {\bf @belapsed} from the {\bf BenchmarkTools.jl}
package to get timings.

We use {\bf lu!} from Julia because then neither factorization will
reallocate the space for the high precision matrix. The excess cost
for the allocation of the low precision matrix in {\bf mplu}
 will show in the timings as will the reduced cost for the factorization.

The problem setup is pretty simple

\begin{code}
julia> using MultiPrecisionArrays

julia> using BenchmarkTools

julia> using MultiPrecisionArrays.Examples

julia> N=4096; G=Gmat(N); A=I - G; x=ones(N); b=A*x;

julia> @belapsed lu!(AC) setup=(AC=copy($A))
1.42840e-01

julia> @belapsed mplu($A)
8.60945e-02

\end{code}

At this point we have timed {\bf lu!} and {\bf mplu}. 
The single precision factorization is a bit more than 
half the cost of the
double precision one.

It is no surprise that the factorization in single precision took roughly half as long as the one in double. In the double-single precision case, iterative refinement is a great example of a time/storage tradeoff. You have to store a low precision copy of , so the storage burden increases by 50\% 
and the factorization time is cut in half. The advantages of IR increase as the dimension increases. IR is less impressive for smaller problems and can even be slower

\begin{code}
julia> N=30; A=I + Gmat(N);

julia> @belapsed mplu($A)
4.19643e-06

julia> @belapsed lu!(AC) setup=(AC=copy($A))
3.70825e-06

\end{code}

Now for the solves. Both {\bf lu} and {\bf mplu}
produce a Julia factorization object and $ \backslash $ works with both.
You have to be a bit careful because MPA and A share storage. So
I will use {\bf lu} instead of {\bf lu!} when factoring $\ma$.

\begin{code}
julia> AF=lu(A); xf = AF\b;

julia> MPAF=mplu!(MPA); xmp=MPAF\b;

julia> luError=norm(xf-x,Inf); MPError=norm(xmp-x,Inf);

julia> println(luError, "  ", MPError)
7.41629e-14  8.88178e-16
\end{code}
So the relative errors are equally good. Now look at the residuals.
\begin{code}
julia> luRes=norm(A*xf-b,Inf)/norm(b,Inf); MPRes=norm(A*xmp-b,Inf)/norm(b,Inf);

julia> println(luRes,"  ",MPRes)
7.40609e-14  1.33243e-15
\end{code}
So, for this well-conditioned problem, IR reduces the factorization cost
by a factor of two and produces results as good as LU on the double precision
matrix. Even so, we should not forget the storage cost of the single
precision copy of $\ma$.

\subsubsection{Harvesting Iteration Statistics: Part I}
\label{sec:stats1}

You can get some iteration statistics by using the
{\bf reporting} keyword argument to the solvers. The easiest way
to do this is with the backslash command. When you use this option you
get a data structure with the solution, the residual norm history, and
the history of the correction norms.

\begin{code}
julia> using MultiPrecisionArrays

julia> using MultiPrecisionArrays.Examples

julia> N=4096; A = I - Gmat(N); x=ones(N); b=A*x;

julia> MPF=mplu(A);

julia> # Use \ with reporting=true

julia> mpout=\(MPF, b; reporting=true);

julia> norm(b-A*mpout.sol, Inf)
2.22045e-16

julia> # Now look at the residual history

julia> mpout.rhist
6-element Vector{Float64}:
 1.00000e+00
 1.64859e-05
 1.02649e-10
 6.35048e-14
 4.44089e-16
 2.22045e-16

julia> # The history of the correction norms is shorter.

julia> mpout.dhist
5-element Vector{Float64}:
 1.00002e+00
 1.62124e-05
 1.02649e-10
 6.35100e-14
 4.45650e-16
\end{code}

One computes a correction after all resduals but the last, when the
iteration terminates on small resdiual norm.

As you can see, IR does well for this problem. The package uses an initial
iterate of $\vx = 0$ and so the initial residual is simply $\vr = \vb$
and the first entry in the residual history is $\| \vb \|_\infty$. The
iteration terminates successfully after four matrix-vector products.

You may wonder why the residual after the first iteration was so much
larger than single precision roundoff. The reason is that the default 
when the low precision is single is to downcast the residual to single before 
the solve (onthefly=false).

One can enable interprecision transfers on the fly and see the difference.

\begin{code}
 
julia> MPF2=mplu(A; onthefly=true);

julia> mpout2=\(MPF2, b; reporting=true);

julia> mpout2.rhist
6-element Vector{Float64}:
 1.00000e+00
 3.10192e-06
 1.04403e-11
 6.13953e-14
 4.44089e-16
 2.22045e-16
\end{code}
There is very little difference in this case.

If we repeat the experiment using half precision as the low precision
the solutions are equally good, but the iteration is slower.

\begin{code}
julia> MPF2=mplu(A; TF=Float16);

julia> # The TF keyword argument lets you make half the low precision.

julia> mpout2 = \(MPF2, b; reporting=true);

julia> norm(A*mpout2.sol - b,Inf)
2.22045e-16

julia> mpout2.rhist
 1.00000e+00
 4.59118e-03
 1.86491e-05
 7.36737e-08
 2.90049e-10
 1.15497e-12
 9.61453e-14
 1.33227e-15
 2.22045e-16
\end{code}

In \S~\ref{sec:exresid} we consider a much more ill-conditioned example.

\subsubsection{Evaluating the residual in extended precision}
\label{sec:exresid}

IR was invented for this case and 
the idea comes from \cite{Wilkinson48}. I am using the notation
from \cite{demmelir,amestoy:2024}. The idea is to evaluate the
residual in a precision TR higher than the working precision. 
One does this by promoting the solution to $TR$ 
and then using the fact that 
if $TR$ is a higher precision than $TW$ and
$\circ$ is any of $+, \times, -, \div$, 
$a \in \calf_R$ and $b \in \calf_W$ that
\begeq
\label{eq:xferrule}
fl_R(a \circ b) = fl_R(a \circ I_W^R(b) ). 
\endeq
This implies that 
\[
\vr = \vb - \ma (I_W^R \vx) = (I_W^R \vb) - (I_W^R \ma) (I_W^R \vx).
\]
After the computation of $\vr$
in precision TR one must downcast it to TW before 
doing the solve for the correction. You must do 
interprecision transfers on the fly in this case. Since the convergence
is based on the norm of the residual in precision TR, you are 
really solving a promoted problem \cite{ctk:irnote}
\[
(I_W^R A) x = I_W^R b
\]
and, by driving the residual to a small value can mitigate ill-conditioning
to some degree.  \pkg allows you to do this with the multiprecision 
factorization you get from {\tt mplu}. The
Krylov-IR solvers can also converge even if the factorization is not
accurate enough to get convergence of IR (see \S~\ref{subsec:irkrylov}).

The classic example is to let TW and TF be single precision and TR be double.
The storage penalty is that you must store two copies of $A$, one for the 
residual computation and the other for the factorization. 
Here is an example with a badly conditioned matrix. You must tell
{\tt mplu} to factor in the working precision and then use the 
{\tt kwargs} in the solver to set TR.

The example in this section is highly ill-conditioned. We base
the example on the integral equation
\begeq
\label{eq:integex}
u - 799.0 \ Gu = f(x) \equiv 1 - \frac{799 x (1-x)}{2}.
\endeq
The solution is $u \equiv 1$. 

\begin{code}
julia> using MultiPrecisionArrays

julia> using MultiPrecisionArrays.Examples

julia> N=4096; alpha=799.0; AD=I - alpha*Gmat(N); 

# conditioning is bad

julia> cond(AD,Inf)
2.35899e+05

# Set up the single precision computation
# and use the right side from the integral equation

julia> h=1.0/(N-1); x=collect(0:h:1); bd = 1.0 .- alpha*x.*(1.0 .- x)*.5;

# Solving in double gives the accuracy you'd expect 

julia> u=AD\bd;

julia> norm(u .- 1.0)
3.16529e-10


# Now move the problem to single precision

julia> A = Float32.(AD); xe=ones(Float32,N); b=Float32.(bd)

julia> # You can see the ill-conditioning

julia> xs = A\b; norm(xs-xe,Inf)
1.37073e-02

julia> # The solution of the promoted problem is better.

julia> xp = Float64.(A)\Float64.(b); norm(xp-xe,Inf)
1.44238e-04

julia> # Make sure TF is what it needs to be for this example

julia> # Set TR in the call to mplu.

julia> AF = mplu(A; TF=Float32, TR=Float64);

julia> # Use the multiprecision array to solve the problem.

julia> mrout = \(AF, b; reporting=true);

julia> # look at the residual history

julia> mrout.rhist
6-element Vector{Float64}:
 9.88750e+01
 6.38567e-03
 1.24560e-06
 9.55874e-08
 8.49478e-08
 8.49478e-08

julia> # Compare the solution to the solution of the promoted problem

julia> norm(mrout.sol - xp,Inf)
5.95629e-08

julia> # That's consistent with theory.

juila> # So the solution is ok?

julia> norm(mrout.sol - xe, Inf)
1.44243e-04

julia> # Yes.
\end{code}

Now I'll try this with TF=half (Float16), the default when
TW = single (Float32). All that one has to do is replace
\begin{code}
AF = mplu(A; TF=Float32, TR=Float64);
\end{code}
with
\begin{code}
AF = mplu(A; TR=Float64);
\end{code}

The iteration fails to converge because $AF$ is not accurate enough.

\begin{code}
julia> AF = mplu(A; TR=Float64);

julia> mout16=\(AF, b; reporting=true);

julia> mout16.rhist
4-element Vector{Float64}:
 9.88750e+01
 3.92451e+00
 3.34292e-01
 2.02204e-01
\end{code}

The advantages of evaluating the residual in extended precision grow
when $A$ is extremely ill-conditioned. Of course, in this case the
factorization in the working precision could be so inaccurate that
IR will fail to converge. One approach to respond to this, as you
might expect, is to use the factorization as a preconditioner and not
a solver \cite{amestoy:2024}. 

\subsection{IR-Krylov}
\label{subsec:irkrylov}

\pkg supports IR-GMRES and IR-BiCGSTAB for TR > TW. You get this to 
work just like in {\tt mplu} by using the keyword argument {\tt TR}.
We will continue with the example in this section and do that. For this 
example the default basis size of $10$ is not enough, so we use 20.

\begin{code}
julia> GF = mpglu(A; TR=Float64, basissize=20);

julia> moutG = \(GF, b; reporting=true);

julia> moutG.rhist
8-element Vector{Float64}:
 9.88750e+01
 7.59075e-03
 1.48842e-05
 2.17281e-07
 8.60429e-08
 7.45077e-08
 7.91866e-08
 7.53089e-08

julia> moutG.dhist
7-element Vector{Float32}:
 1.03306e+00
 3.23734e-02
 2.83073e-04
 8.92629e-06
 1.55432e-07
 6.05685e-08
 5.96633e-08


# We got the solution to the promoted problem...

julia> xp = Float64.(A)\b;

julia> norm(xp-moutG.sol, Inf)
5.95825e-08

# IR-BiCGSTAB would take fewer iterations than IR-GMRES had we used
# the default basis size because there's no storage
# issue. But remember that BiCGSTAB has a higher cost per linear iteration.

julia> BF = mpblu(A; TR=Float64);

julia> moutB = \(BF, b; reporting=true);

julia> moutB.rhist
5-element Vector{Float64}:
 9.88750e+01
 1.86437e-07
 7.53089e-08
 7.53089e-08
 7.53089e-08

julia> moutB.dhist
4-element Vector{Float32}:
 1.00227e+00
 2.05457e-06
 5.95821e-08
 5.95821e-08

julia> norm(xp - moutB.sol, Inf)
5.95825e-08
\end{code}

\subsection{Options and data structures for mplu}
\label{subsec:options}

Here is the source for {\tt mplu}.

\begin{code}
"""
function mplu(A::AbstractArray{TW,2}; TF=nothing, TR=nothing, 
       residterm=residtermdefault, onthefly=nothing) where TW <: Real
#
# If the high precision matrix is single, the low precision must be half
# unless you're planning on using a high-precision residual where TR > TW
# and also factoring in the working precision, so TW == TF.
# 
# 
(TR == nothing) && (TR = TW)
TFdef = Float32
(TW == Float32) && (TFdef = Float16)
(TF == nothing) && (TF = TFdef)
#
# Unless you tell me otherwise, onthefly is true if low precision is half
# and false if low precision is single.
#
(onthefly == nothing ) && (onthefly = (TF==Float16)) 
#
# IF TF = TW then something funny is happening with the residual precision.
#
(TF == TW) && (onthefly=true)
#
# Build the multiprecision array MPA
#   
MPA=MPArray(A; TF=TF, TR=TR, onthefly=onthefly)
#   
# Factor the low precision copy to get the factorization object MPF
#   
MPF=mplu!(MPA; residterm=residterm) 
return MPF
end 
"""
\end{code}

The function {\tt mplu} has four keyword arguments. The easy one to 
understand is {\tt TF} which is the precision of the factorization. Julia 
has support for single ({\tt Float32}) and half ({\tt Float16})
precisions. If you set {\tt TF=Float16} then low precision will be half. 
Don't do that unless you know what you're doing. Using half precision 
is a fast way to get incorrect results. Look at \S~\ref{sec:half}
for a bit more bad news.

{\tt TR} is the residual precision. The {\tt MPArray} structure preallocates
storage for the residual and a local copy of the solution in precision 
{\tt TR}. For most applications (but see \S~\ref{sec:exresid}) 
{\tt TR = TW}.

The final keyword arguments are {\tt onthefly} and {\tt residterm}.
The {\tt onthefly} keyword controls 
how the triangular solvers from the factorization work. When you solve

$$
\ml \mU \vd = \vr
$$

the LU factors are in low precision and the residual $r$ 
is in high precision. If you let Julia and LAPACK figure out what to do, 
then the solves will be done in high precision and
the entries in the LU factors will be converted to high precision with 
each binary operation. The output $d$ will be in high precision. 
This is called interprecision transfer on-the-fly
and {\tt onthefly = true} will tell the solvers to do it that way. 
You have $N^2$ interprecision transfers with each solve and, 
as we will see, that can have a non-trivial cost.

When low precision is Float32 and {\tt TR = TW = Float64} 
then the default is {\tt onthefly = false}. 
I am rethinking this in
view of the experiments in \cite{ctk:irnote} and make 
{\tt onthefly = true} the default pretty soon.

When you set {\tt onthefly = false} the solve
converts $r$ to low precision, does the solve entirely in 
low precision, and then promotes $d$ to high precision. 
This is {\bf in-place} interprecision transfer. 
You need to 
be careful to avoid overflow and, more importantly, underflow when 
you do that and we scale $r$ to be a unit vector before conversion 
to low precision and reverse the scaling when we promote $d$. 
We take care of this for you.

The {\tt residterm} keyword controls the termination options. See
\S~\ref{subsec:termination} for the details on that. I will be 
making changes to the way \pkg controls termination frequently.

{\tt mplu} calls the constructor for the multiprecision array and then 
factors the low precision matrix. In some cases, such as nonlinear solvers, 
you will want to separate the constructor and the factorization. I will
provide code for this in a future release. I do not export the constructor
and if you mess with it, remember that {\tt mplu!} 
overwrites the low precision copy of A with the factors. The factorization
object is different from the multiprecision array, even though they
share storage. Be careful with this.

\subsection{Memory Allocations: I}
\label{subsec:alloc1}

The memory footprint of a multiprecision array is dominated by
the high precision array and the low precision copy. The allocations of
\begin{code}
AF1=lu(A)
\end{code}
and
\begin{code}
AF2=mplu(A)
\end{code}
are very different. Typically {\tt lu} makes a high precision copy of $\ma$ and
factors that with {\tt lu!}. {\tt mplu} on the other hand, uses $\ma$
as the high precision matrix in the multiprecision array structure and
the makes a low precision copy to send to {\tt lu!}. Hence {\tt mplu}
has half the allocation burden of {\tt lu}.

That is, of course misleading. The most memory-efficient 
way to apply {\tt lu} is to
overwrite $\ma$ with the factorization using
\begin{code}
AF1=lu!(A).
\end{code}
The analog of this approach with a multiprecision array would be to
first build an {\tt MPArray} structure with a call to the constructor
\begin{code}
MPA = MPArray(A)
\end{code}
which makes $\ma$ the high precision matrix and also makes a low
precision copy. This is the stage where the extra memory is allocated
for the the low precision copy. One follows that with the factorization
of the low precision matrix to construct the factorization object.
\begin{code}
MPF = mplu!(MPA).
\end{code}
The function {\tt mplu} simply applies {\tt MPArray} and follows that with
{\tt mplu!}. 

I do not export the constructor {\tt MPArray}. You should use {\tt mplu}
and not import the constructor.

Once you have used {\tt mplu} to make a multiprecision factorization, you 
can reuse that storage for a different matrix as long as the size
and the precision are the same. For example, suppose
\begin{code}
MPF = mplu(A)
\end{code}
is a multiprecision factorization of $\ma$. If you want to factor $\mb$
and reuse the memory, then 
\begin{code}
MPF = mplu!(MPF,B)
\end{code}
will do the job.

\section{Half Precision} 
\label{sec:half}

There are two half precision (16 bit) formats. Julia has native
support for IEEE 16 bit floats (Float16). A second format
(BFloat16) has a larger exponent field and a smaller 
significand (mantissa), thereby
trading precision for range. In fact, the exponent field in BFloat is
the same size (8 bits) as that for single precision (Float32). The 
significand, however, is only 8 bits. Which is less than that for 
Float16 (11 bits) and single (24 bits). The size of the significand
means that you can get in real trouble with half precision in either
format. 

At this point Julia has no native support for BFloat16. 
Progress is being made and we expect to support BFloat16
in the future.

Doing the factorization in half precision (Float16 or BFloat16) 
will not speed up the solver, in fact it will make the 
solver slower. The reason for this is that LAPACK and the BLAS 
do not ({\bf YET} \cite{newblas}) 
support half precision, so all the clever stuff in
there is missing. We provide a half precision LU factorization 
for Float16
{\bf /src/Factorizations/hlu!.jl} that is better than nothing. 
It's a hack of Julia's  {\tt generic\_lu!} with threading and a couple
compiler directives. Even so, it's 3 -- 7 times
{\bf slower} than a double precision LU. Half precision support is coming 
\cite{newblas} and 
Julia and Apple support Float16 in hardware. Apple also has hardware
support for BFloat16. However,
for now, at least for
desktop computing, half precision is for
research in iterative refinement, not applications.

Here's a table 
that illustrates the point. In the table we compare timings
for LAPACK's LU to the LU 
we compute with {\tt hlu!.jl}. The matrix is $\mi - 800.0 \mg$.

\begin{table}[h!]
\label{tab:halftime}
\caption{Half precision is slow: LU timings}
\centerline{
\begin{tabular}{lllll} 
        N &   Double &   Single &     Half &    Ratio \\ 
\hline 
512 & 1.37e-03 & 6.32e-04 & 3.84e-03 & 2.82e+00   \\ 
1024 & 4.97e-03 & 2.82e-03 & 3.51e-02 & 7.07e+00   \\ 
2048 & 2.67e-02 & 1.49e-02 & 1.91e-01 & 7.15e+00   \\ 
4096 & 1.64e-01 & 8.32e-02 & 7.58e-01 & 4.63e+00   \\ 
8192 & 1.25e+00 & 7.18e-01 & 5.67e+00 & 4.53e+00   \\ 
\hline 
\end{tabular} 
}
\end{table}

The columns of the table are the dimension of the problem, timings
for double, single, and half precision, and the ratio of the half
precision timings to double. The timings came from Julia 1.10-beta2
running on an Apple M2 Pro with 8 performance cores.

Half precision is also difficult to use properly. 
The low precision can make iterative refinement fail because the 
half precision factorization can have
a large error. Here is an example to illustrate this point. 
The matrix here is modestly ill-conditioned and you can see that 
in the error from a direct solve in double precision.

\begin{code}
julia> A=I - 800.0*G;

julia> x=ones(N);

julia> b=A*x;

julia> xd=A\b;

julia> norm(b-A*xd,Inf)
6.96332e-13

julia> norm(xd-x,Inf)
2.30371e-12
\end{code}

Now, if we downcast things to half precision, nothing good happens.

\begin{code}
julia> AH=Float16.(A);

julia> AHF=hlu!(AH);

julia> z=AHF\b;

julia> norm(b-A*z,Inf)
6.25650e-01

julia> norm(z-xd,Inf)
2.34975e-01
\end{code}

So you get very poor, but unsurprising, results. While MultiPrecisionArrays.jl supports half precision and I use it all the time, it is not something you would use in your own work without looking at the literature and making certain you are prepared for strange results. Getting good results consistently from half precision is an active research area.

So, it should not be a surprise that IR also struggles with half precision.
We will illustrate this with one simple example. In this example high
precision will be single and low will be half. Using {\bf MPArray} with
a single precision matrix will automatically make the low precision matrix
half precision. In this example we use the keyword argument 
```onthefly``` to toggle between MPS and LPS.
\begin{code}
julia> N=4096; G=800.0*Gmat(N); A=I - Float32.(G);

julia> x=ones(Float32,N); b=A*x;

julia> MPF=mplu(A; onthefly=false);

julia> y=MPF\b;

julia> norm(b - A*y,Inf)
1.05272e+02
\end{code}
So, IR completely failed for this example. We will show how to extract
the details of the iteration in a later section.

It is also worthwhile to see if doing the triangular solves on-the-fly 
(MPS) helps. 

\begin{code}
julia> MPF2=mplu(A; onthefly=true);

julia> z=MPF2\b;

julia> norm(b-A*z,Inf)
1.28174e-03
\end{code}
So, MPS is better in the half precision case. Moreover, it is also less
costly thanks to the limited support for half precision computing.
For that reason, MPS is the default when high precision is single.

However, on-the-fly solves are not enough to get good results and IR
still terminates before converging to the correct result.

\section{Using the Low Precision Factorization as a Preconditioner}
\label{sec:precond}

In this section we present some options if IR fails to converge. This 
is very unlikely if high precision is double and low precision is single.
If low precision is half, the methods in this section might save you.

The idea is simple. Even if 
\[
\mm_{IR} = \mi - {\hat \mU}^{-1} {\hat \ml}^{-1} \ma
\]
has norm larger than one, it could still be the case that
\[
{\hat \mU}^{-1} {\hat \ml}^{-1} \ma
\]
is well conditioned and that 
\begeq
\label{eq:pdef}
\mP={\hat \mU}^{-1} {\hat \ml}^{-1}
\endeq
could be a useful preconditioner for a Krylov method. 

\subsection{Direct Preconditioning}
\label{subsec:preconddirect}

The obvious way to use $\mP$ is simply to precondition the equation
$\ma \vx = \vp$. In this case we prefer right preconditioning where we 
solve
\[
\ma \mP \vz = \vb
\]
and then set $\vx = \mP \vx$. This is different from all IR methods
we discuss in this paper and one may lose some accuracy by avoiding the
IR loop. 

\subsection{Krylov-IR}
\label{subsec:kyrir}

Krylov-IR methods solve the correction equation with a preconditioned
Krylov iteration using the low precision solve as the preconditioner. 
Currently \pkg supports GMRES \cite{gmres} and BiCGSTAB
\cite{bicgstab}.

\subsection{GMRES-IR}
\label{subsec:gmresir}

GMRES-IR \cite{CarsonHigham1,CarsonHigham} solves the correction equation
with a preconditioned GMRES \cite{gmres} iteration. One way to think of this
is that the solve in the IR loop is an approximate solver for the
correction equation
\[
\ma \vd = \vr
\]
where one replaces $\ma$ with the low precision factors
$\ml \mU$. In GMRES-IR one solves the correction
equation with a left-preconditioned GMRES iteration using 
$\mP$ as
the preconditioner. The preconditioned equation is
\[
\mP \ma \vd = \mP \vr.
\]
The reason for using left preconditioning is that one is not
interested in a small residual for the correction equation, but in
capturing $\vd$ as well as possible. The IR loop is the part of the solve
that seeks a small residual norm.

GMRES-IR will not be as efficient as IR because each iteration is itself
an GMRES iteration and application of the preconditioned matrix-vector
product has the same cost (solve + high precision matrix vector product)
as a single IR iteration. However, if low precision is half, this approach
can recover the residual norm one would get from a successful IR iteration.

There is also a storage problem. One should allocate storage for the Krylov
basis vectors and other vectors that GMRES needs internally. We do that
in the factorization phase. So the structure {\bf MPGEFact} has the 
factorization of the low precision matrix, the residual, the Krylov
basis and some other vectors needed in the solve. The Julia function
{\bf mpglu} constructs the data structure and factors the low precision
copy of the matrix. The output, like that of {\bf mplu} is a factorization
object that you can use with backslash.

Here is a well conditioned example. Both IR and GMRES-IR perform well, with
GMRES-IR taking significantly more time. In these examples high precision
is single and low precision is half. 

\begin{code}
julia> using MultiPrecisionArrays

julia> using MultiPrecisionArrays.Examples

julia> using BenchmarkTools

julia> N=4069; AD= I - Gmat(N); A=Float32.(AD); x=ones(Float32,N); b=A*x;

julia> # build two MPArrays and factor them for IR or GMRES-IR

julia> MPF=mplu(A); MPF2=mpglu(A);

julia> z=MPF\b; y=MPF2\b; println(norm(z-x,Inf),"  ",norm(y-x,Inf))
5.9604645e-7  4.7683716e-7

julia> # and the relative residuals look good, too

julia> println(norm(b-A*z,Inf)/norm(b,Inf),"  ",norm(b-A*y,Inf)/norm(b,Inf))
4.768957e-7  3.5767178e-7

julia> @btime $MPF\$b;
  13.582 ms (4 allocations: 24.33 KiB)

julia> @btime $MPF2\$b;
  40.028 ms (183 allocations: 90.55 KiB)
\end{code}

If you dig into the iteration statistics (more on that later) you
will see that the GMRES-IR iteration took almost exactly four times
as many solves and residual computations as the simple IR solve.

We will repeat this experiment on the ill-conditioned example. In this
example, as we saw earlier, IR fails to converge.

\begin{code}
julia> N=4069; AD= I - 800.0*Gmat(N); A=Float32.(AD); x=ones(Float32,N); b=A*x;

julia> MPF=mplu(A); MPF2=mpglu(A);

julia> z=MPF\b; y=MPF2\b; println(norm(z-x,Inf),"  ",norm(y-x,Inf))
0.2875508  0.0044728518

julia> println(norm(b-A*z,Inf)/norm(b,Inf),"  ",norm(b-A*y,Inf)/norm(b,Inf))
0.0012593127  1.4025759e-5

\end{code}

So, the relative error and relative residual norms for GMRES-IR
are much smaller than for IR.

\subsection{Memory Allocations: II}
\label{subsec:alloc2}

Much of the discussion from \S~\ref{subsec:alloc1} remains valid for
the {\tt MGPArray} structure and the associated factorization
structure {\tt MPGEFact}. The only difference that matters is that
{\tt MGPArray} contains the Krylov basis and a few other vectors
that GMRES needs, so the allocation burden is a little worse. 

The allocation burden is less for BiCGSTAB-IR and the {\tt MPBArray}
structure.  

That aside, {\tt mpglu!} and {\tt mpblu!} 
work the same way that {\tt mplu!} does
when factoring or updating a {\tt MGPArray}.

The functions {\tt mpglu} and {\tt mpblu} combine the allocations and
the factorization.

\subsection{Harvesting Iteration Statistics: Part 2}
\label{sec:stats2}

The output for GMRES-IR contains the residual history and a vector with the
number of Krylov iterations for each IR step. The next
example illustrates that.

\begin{code}
julia> MPGF=mpglu(A);

julia> moutg=\(MPGF, b; reporting=true);

julia> norm(A*moutg.sol-b, Inf)
1.44329e-15

julia> moutg.rhist
3-element Vector{Float64}:
 9.99878e-01
 7.48290e-14
 1.44329e-15

julia> moutg.khist
2-element Vector{Int64}:
 4
 4
\end{code}

While only two IR iterations are needed for convergence, the Krylov history
shows that each of those IR iterations needed four GMRES iterations. Each
of those GMRES iterations requires a matrix-vector product and a low-precision
on-the-fly linear solve. So GMRES-IR is more costly and, as pointed out
in \cite{CarsonHigham1,CarsonHigham} is most useful with IR does not converge
on its own.

We will demonstrate this with one last example. In this example high precision
is single and low precision is half. As you will see, this example is 
very ill-conditioned.

\begin{code}

julia> N=8102; AD = I - 799.0*Gmat(N); A=Float32.(AD); x=ones(Float32,N); b=A*x;

julia> cond(A, Inf)
2.34824e+05

julia> MPFH=mplu(A);

julia> mpouth=\(MPFH, b; reporting=true);

julia> # The iteration fails. 

julia> mpouth.rhist
4-element Vector{Float64}:
 9.88752e+01
 9.49071e+00
 2.37554e+00
 4.80087e+00

julia> # Try again with GMRES-IR and mpglu

julia> MPGH=mpglu(A);

julia> mpoutg=\(MPGH, b; reporting=true);

julia> mpoutg.rhist
4-element Vector{Float32}:
 9.88752e+01
 1.86920e-03
 1.29700e-03
 2.46429e-03

julia> mpoutg.khist
3-element Vector{Int64}:
 10
 10
 10

\end{code}

So GMRES-IR does much better. Note that we are taking ten GMRES iterations
for each IR step. Ten is the default. To increase this set the
keyword argument {\bf basissize}.

\section{Details}
\label{sec:details}

In this section we discuss a few details that are important for 
understanding IR, but less important for simply using \pkg.

\subsection{Terminating the while loop}
\label{subsec:termination}

There are many parameters in the termination criteria and the
defaults in \pkg are reasonable for most applications. These
defaults differ slightly from the recommendations in \cite{higham97}
and we explain that later in this section.

One can terminate the iteration when the relative residual is small
\ie when
\begeq
\label{eq:residterm}
\| \vr \| < C_r u_w \| \vb \|
\endeq
or when the normwise backward error is small, \ie
\begeq
\label{eq:term}
\| \vr \| < C_e u_w  (\| \vb \| + \| \ma \| \| \vx \|).
\endeq
You make the choice
between \eqnok{residterm} and \eqnok{term} when you compute the
multiprecision factorization. 

The prefactors $C_r$ and $C_e$ are algorithmic parameters
$C_r = 1.0$ is the default in \pkg.
The recommendation in 
\cite{higham97} is that $C_e = 1.0$, which is the default in \pkg.

If {\tt TR} is wider than {\tt TW}, then I assume that you are 
solving a highly ill-conditioned problem and terminate when the
norms of the corrections stagnate
\begeq
\label{eq:wilkterm}
\| \vd_{new} \| \ge R_{max} \| \vd_{old} \|.
\endeq
This is a bit different from the approach in \cite{demmelir} and
I may change to that way in the future. Using \eqnok{wilkterm}
is simpler, but might take one more iteration to drive the error
to $u_w$, as the theory predicts if IR converges.

\subsubsection{Residual or Backward Error}

Using the normwise backward
error is more costly because the computation of $\| \ma \|$ is
$N^2$ work and is more expensive that a few IR iterations. 

In \pkg we do this when we compute $\| \ma \|$ during
the multiprecision factorization. After we copy $\ma$ to the factorization
precision, we can take the norm of the low-precision copy before we factor
it. Specifically, we
compute $\| \ma \|_1$ in the factorization precision if 
{\tt TF = Float32} or {\tt TF = Float64}. If {\tt TF} is 
half precision, which is still slow in our desktop environment, we
compute $\| \ma \|$ in the working precision {\tt TW}. We control this
via the {\tt residterm} {\tt kwarg} for the factorizations. 

The user can change from using \eqnok{residterm} to \eqnok{term} with 
the keword parameter {\tt residterm} in {\tt mplu}. 
The default value is {\tt true}.
Setting it to {\tt false} enables termination on small normwise
backward error. For example
\begin{code}
MPF = mplu(A; residterm=false)
\end{code}
tells the solver to use \eqnok{term}.

\subsubsection{Protection Against Stagnation}

The problem with either of \eqnok{residterm} or \eqnok{term} is
that IR can stagnate, especially for ill-conditioned problems, before
the termination criterion is attained. We detect stagnation by looking
for a unacceptable decrease (or increase) in the residual norm. So we will
terminate the iteration if
\begeq
\label{eq:stag}
\| \vr_{new} \| \ge R_{max} \| \vr_{old} \|
\endeq
even if the small residual condition is not satisfied.
The recommendation in \cite{higham97} is $R_{max} = .5$ and that is
the default in \pkg.

The paper \cite{higham97} also recommends that on put a limit 
{\tt litmax = 5} on
the number of IR iterations. We use a higher
value of {\tt litmax = 10} in \pkg.

In this paper we count iterations as residual computations. This means
that the minimum number of iterations will be two. Since we begin
with $\vx = 0$ and $\vr = \vb$, the first iteration (iteration $0$) computes
$\vd = \mU^{-1} \ml^{-1} \vb$ and then $\vx \leftarrow \vx + \vd$, so
the first iteration is the output of a low precision solve.  We will
need at most one more iteration to test for a meaningful residual reduction.
Iterations after iteration $0$ cost more because they do a meaningful 
residual computation.

\subsubsection{Changing the parameters}

The parameters and defaults are
\begin{itemize}
\item $C_r = 1.0$
\item $C_e = 1.0$
\item $R_{max} = .5$
\item $litmax = 10$
\end{itemize}

We store the parameters in a TERM structure, which define in the main file
for \pkg.
\begin{code}
struct TERM
       Cr::Real
       Ce::Real
       Rmax::Real
       litmax::Int
end
\end{code}
We create one {\tt TERM} structure for the defaults
{\tt term\_parms\_default}. The solvers take a TERM structure 
as a kwarg, with {\tt term\_parms\_default} as the default.

To change the parameters use the {\tt update\_parms} function. This
function makes a TERM structure for you to pass to the solvers.
\begin{code}
function update_parms(; Cr=Cr_default, Ce=Ce_default,
      Rmax=Rmax_default, litmax=litmax_default)
term_parms=TERM(Cr, Ce, Rmax, litmax)
return term_parms
end
\end{code}

Here is an example with the ill-conditioned problem we've been using
in other examples. In this example we decrease Rmax and see that the
solve takes fewer iterations with no change in the residual quality.
\begin{code}
julia> using MultiPrecisionArrays

julia> using MultiPrecisionArrays.Examples

julia> N=512; A=I - 799.0*Gmat(N); AF=mplu(A);

julia> b=ones(N);

julia> mout=\(AF,b; reporting=true);

# Termination on a residual norm increase.

julia> mout.rhist
6-element Vector{Float64}:
 1.00000e+00
 4.39096e-03
 2.85170e-07
 4.30167e-11
 6.05982e-12
 6.34648e-12

julia> term_parms=update_parms(;Rmax=.1);

julia> mout2=\(AF,b; reporting=true, term_parms=term_parms);

# Termination on minimal decrease in the residual norm. 

julia> mout2.rhist
5-element Vector{Float64}:
 1.00000e+00
 4.39096e-03
 2.85170e-07
 4.30167e-11
 6.05982e-12

\end{code}

\subsection{Is $O(N^2)$ really negligible}
\label{subsec:nsquared}

In this section {\tt TR = TW = Float64} and 
{\tt TF = TS = Float32}, which means that the iterprecision transfers
in the triangular solvers are done in-place. We terminate on small residuals.

The premise behind IR is that reducing the $O(N^3)$ cost of the
factorization will make the solve faster because everything else
is $O(N^2)$ work. It's worth looking into this.

We will use the old-fashioned definition of a FLOP as an add, a multiply
and a bit of address computation. So we have $N^2$ flops for any of
(1) matrix-vector multiply $\ma \vx$, 
(2) the two triangular solves with the LU factors
$(\ml \mU)^{-1} \vb$, and 
(3) computation of the $\ell^1$ or $\ell^\infty$ matrix operator norms
$\| \ma \|_{1,\infty}$.

A linear solve with an LU factorization and the standard triangular
solve has a cost of $(N^3/3) + N^2$ TR-FLOPS. The factorization for IR
has a cost of $N^3/3$ TF-FLOPS or $N^3/6$ TR-FLOPS.

A single IR iteration costs a matrix-vector product in precision TR
and a triangular solve in precision TF for a total of
$3 N^2/2$ TR-FLOPS. Hence a linear solve with IR that needs $n_I$ iterations
costs
\[
\frac{N^3}{6} + 3 n_I N^2/2
\]
TR-FLOPS if one terminates on small residuals and an extra $N^2$ TR-FLOPS
if one computes the norm of $\ma$ in precision TR.

IR will clearly be better for large values of $N$. How large is that?
Table~\ref{tab:solvecomp} and Table~\ref{tab:solvecomp2}
compare the time for factorization 
nd solve
using $lu!$ and $ldiv$ (cols 2-4) with the equivalent multiprecision
commands $mplu$ and $\backslash$. The final column is the time for
computing $\| \ma \|_1$, which would be needed to terminate on
small normwise backward errors.

\clearpage

The operator is {\tt A = I - Gmat(N)}. We tabulate

\begin{itemize}
\item LU: time for {\tt AF=lu!(A)}
\item TS: time for {\tt ldiv!(AF,b)}
\item TOTL = LU+TS
\item MPLU: time for {\tt MPF=mplu(A)}
\item MPS: time for {\tt MPF\textbackslash  b}
\item TOT: MPLU+MPS
\item OPNORM: Cost for $\| A \|_1$, which one needs to terminate on small normwise backward error.
\end{itemize}

The message from the tables is that the triangular solves are more costly
than operation counts might indicate. One reason for this is that the
LU factorization exploits multi-core computing better than a triangular
solve. It is also interesting to see how the choice of BLAS affects the
relative costs of the factorization and the solve. Also notice
how the computation of $\| \ma \|$ is more costly than the
two triangular solves for {\tt ldiv!(AF,b)}.

For both cases, multiprecision arrays perform better when $N \ge 2048$
with the difference becoming larger as $N$ increases.

\begin{table}[h!]
\caption{\label{tab:solvecomp}$\alpha=1$, openBLAS}
\centerline{
\begin{tabular}{llllllll} 
        N &       LU &       TS &      TOTL &     MPLU &      MPS &      TOT &   OPNORM \\ 
\hline 
512 & 9.8e-04 & 6.4e-05 & 1.0e-03 & 7.8e-04 &  2.9e-04 & 1.1e-03 &  1.7e-04   \\ 
1024 & 3.7e-03 & 2.6e-04 & 4.0e-03 & 3.1e-03 &  8.3e-04 & 4.0e-03 &  7.8e-04   \\ 
2048 & 2.1e-02 & 1.5e-03 & 2.3e-02 & 1.4e-02 &  5.4e-03 & 1.9e-02 &  3.4e-03   \\ 
4096 & 1.5e-01 & 5.4e-03 & 1.5e-01 & 8.8e-02 &  2.0e-02 & 1.1e-01 &  1.4e-02   \\ 
8192 & 1.1e+00 & 2.1e-02 & 1.1e+00 & 6.0e-01 &  7.3e-02 & 6.7e-01 &  5.8e-02   \\ 
\hline 
\end{tabular} 
}
\end{table}

\begin{table}[h!]
\caption{\label{tab:solvecomp2}$\alpha=1$, AppleAcclerateBLAS}
\centerline{
\begin{tabular}{llllllll} 
        N &       LU &       TS &      TOTL &     MPLU &      MPS &      TOT &   OPNORM \\ 
\hline 
512 & 7.9e-04 & 1.1e-04 & 9.0e-04 & 4.7e-04 &  2.9e-04 & 7.6e-04 &  1.7e-04   \\ 
1024 & 3.4e-03 & 4.9e-04 & 3.9e-03 & 2.4e-03 &  1.2e-03 & 3.6e-03 &  7.8e-04   \\ 
2048 & 1.9e-02 & 2.4e-03 & 2.1e-02 & 1.1e-02 &  1.1e-02 & 2.2e-02 &  3.4e-03   \\ 
4096 & 1.4e-01 & 1.2e-02 & 1.5e-01 & 6.1e-02 &  5.6e-02 & 1.2e-01 &  1.4e-02   \\ 
8192 & 1.2e+00 & 5.5e-02 & 1.3e+00 & 5.8e-01 &  2.2e-01 & 8.0e-01 &  5.7e-02   \\ 
\hline 
\end{tabular} 
}
\end{table}

\subsection{Interprecision Transfers: Part I}
\label{subsec:interprec1}

The meaning of $\vd = \mU^{-1} \ml^{-1} \vr$ is more subtle. The problem
is that the factors $\mU$ and $\ml$ are store in low precision and $\vr$
is a high precision vector. LAPACK will convert $\ml$ and $\mU$ to the
higher precision ``on the fly'' with each mixed precision binary operation
at a cost of $O(N^2)$ interprecision transfers. 
The best way to understand this is from \eqnok{xferrule}, which implies that
if $\ml$ and $\mU$ are in precision $TF$ and $\vb$ is in precision $TW$, then
\[
(\ml \mU)^{-1} \vb = ( (I_F^W \ml) (I_F^2 \mu) )^{-1} \vb.
\]
As we will see this interprecision transfer can have a meaningful cost
even though the factorization will dominate with $O(N^3)$ work.

One can eliminate the the $O(N^2)$ interprecision transfer
cost by copying $\vr$ into low precision, doing
the triangular solves in low precision, and then mapping the result
into high precision. The two approaches are not the same.
To see this we
$\vx_c$ denote the current iterate and $\vx_+$ the new iterate.

If one does
the solves on the fly then the IR iteration
\[
\begin{array}{ll}
\vx_+ & = \vx_c + \vd
= \vx_c + {\hat \mU}^{-1} {\hat \ml}^{-1} \vr\\
\\
& = \vx_c + {\hat \mU}^{-1} {\hat \ml}^{-1} (\vb - \ma \vx_c) \\
\\
& = (\mi - {\hat \mU}^{-1} {\hat \ml}^{-1} \ma ) \vx_c
+ {\hat \mU}^{-1} {\hat \ml}^{-1} \vb
\end{array}
\]
is a linear stationary iterative method. Hence on the fly IR
will converge if the spectral radius of the iteration matrix
\[
\mm_{IR} = \mi - {\hat \mU}^{-1} {\hat \ml}^{-1} \ma
\]
is less than one. We will refer to the on the fly approach as
mixed precision solves (MPS) when we report computational
results in \S~\ref{subsec:precision2}. In the terminology of
\cite{amestoy:2024}, the MPS approach sets the solver precision $TS$
to the precision of the residual computation $TR$, which is $TW$ for
now.

If one does the triangular solves in low precision,
one must first take care to scale $\vr$ to avoid underflow, so one solves
\begeq
\label{eq:vdlow}
(\ml \mU) \vd_F = I_W^F(\vr/\| \vr \|)
\endeq
in low precision and then promotes $\vd$ to high precision and reverses
the scaling to obtain
\begeq
\label{eq:promotevd}
\vd = \| \vr \| I_F^W (\vd_F).
\endeq

We will refer to this approach as low precision solves (LPS) when
we report computational results in \S~\ref{subsec:precision2}.
In practice, if low precision is single, the quality of the results
is as good as one would get with MPS and the solve phase is somewhat
faster. Making the connection to \cite{amestoy:2024} again, the LPS
approach sets $TS = TF$.

We express this choice in \pkg with the {\tt onthefly} keyword 
argument which tells the solvers to do the triangular solves on the
fly ({\tt true}) or not ({\tt false}). The default is {\tt true} unless
$TW$ is double and $TF$ is single.

\subsection{Convergence Theory}
\label{sec:convergence}

\subsubsection{Estimates for $\| \mm_{IR} \|$}
\label{subsec:normest1}

We will estimate the norm of $\mm_{IR}$ to see how the factorization
precision affects the convergence. First write
\begeq
\label{eq:mexp}
\mm_{IR} = \mi - {\hat \mU}^{-1} {\hat \ml}^{-1} \ma
= {\hat \mU}^{-1} {\hat \ml}^{-1} ( {\hat \ml} {\hat \mU} - \ma ).
\endeq
We split $\Delta \ma = ( {\hat \mU} {\hat \ml} - \ma )$ to separate
the rounding error from the backward error in the low precision
factorization
\[
\Delta \ma = ({\hat \mU} {\hat \ml} - I_W^L \ma) +
(I_W^L \ma - \ma).
\]
The last term can be estimated easily
\begeq
\label{eq:roundmm}
\| I_W^F \ma - \ma \| \le u_F \| \ma \|.
\endeq

To estimate the first term we look at the component-wise
backward error \cite{higham}. If $3Nu_F < 1$ then
\begeq
\label{eq:comperr}
| {\hat \ml} {\hat \mU} - I_W^L \ma | \le \gamma_{3N} (u_F)
| {\hat \ml} | | {\hat \mU} |.
\endeq
In \eqnok{comperr} $| \mb |$ is the matrix with entries
the absolute values of those in $\mb$ and
\[
\gamma_k(u) = \frac{ku}{1 - ku}.
\]
We can combine \eqnok{roundmm} and \eqnok{comperr} to get 
\begeq
\label{eq:bottomline1}
\begin{array}{ll}
\| \mm_{IR} \| & \le u_F \| \ma \| + \gamma_k(u_F)  
\| {\hat \mU}^{-1} {\hat \ml}^{-1} \|
\| {\hat \ml} \| \| {\hat \mU} \|\\
\\
& = u_F \| \ma \| + \gamma_k(u_F) 
\kappa({\hat \ml})\kappa(\hat \mU).
\end{array}
\endeq
The standard estimate in textbooks 
for $\| {\hat \ml} \| \| {\hat \mU} \|$
uses very pessimistic
(and unrealistic) worst case bounds on the right side of \eqnok{comperr}.
In cases where the conditioning of the factors is harmless, 
the estimate in \eqnok{bottomline1} suggests that IR should converge
well if low precision is single.

We will use the probabilistic bounds from \cite{HighamMary} to explore
this in more detail. Roughly
speaking, with high probability for
desktop sized $N \le 10^{10}$ problems we obtain
\begeq
\label{eq:highammary}
| {\hat \ml} {\hat \mU} - I_W^L \ma | \le (13 u_F \sqrt{N} + O(u_F^2))
\| {\hat \ml} \| \| {\hat \mU} \|.
\endeq
If we neglect the $O(u_F^2)$ term in \eqnok{highammary}, our estimate
for $\mm$ becomes
\begeq
\label{eq:bottomline2}
\begin{array}{ll}
\| \mm \| & \le 
u_F (\| \ma \| + \| {\hat \mU}^{-1} {\hat \ml}^{-1} \|
13 \sqrt{N} \| {\hat \ml} \| \| {\hat \mU} \|) \\
\\
&\le u_l ( \| \ma \| + 13 \sqrt{N} \kappa({\hat \ml})\kappa(\hat \mU) ).
\end{array}
\endeq
So, $\| \mm \| < 1$ if
\[
\| \ma \| + 13 \sqrt{N} \kappa({\hat \ml})\kappa(\hat \mU) < u_F^{-1}.
\]
For example if we assume that $\| \ma \| =O(1)$, low precision is single
($u_F = u_s \approx 1.2 \times 10^{-7}$),
and we make a fairly pessimistic assumption about the
conditioning of the low precision factors,
\[
\kappa({\hat \ml})\kappa(\hat \mU) \le \sqrt{N},
\] then $\| \mm \| < 1$ if
\begeq
\label{eq:nbounds}
N < u_s^{-1}/14 \approx 6 \times 10^5
\endeq
which is the case for most desktop sized problems. However, if low
precision is half, then \eqnok{nbounds} becomes with
$u_F = u_h \approx 9.8 \times 10^{-3}$
\begeq
\label{eq:nboundh}
N < u_F^{-1}/14 \approx 73.
\endeq
This is an indication that there are serious risks in using half
precision if the conditioning of the low precision factors
increases with $N$, which could be the case if the $\ma$ is a
discretization of a boundary value problem.

\subsubsection{Limiting Behavior of IR}
\label{subsec:limit}

In exact arithmetic one would get a reduction in the error with
each iteration of a factor of 
$\rho(\mm_{IR}) \le \| \mm_{IR} \|$. However, when one
accounts for the errors in the residual computation, we will see
how and when the iteration can stagnate. Our analysis will be a
simplified version of the one from \cite{CarsonHigham} and we will
neglect many of the details.

In this section we will consider dense matrices with solves with
MPS, so the solves with the low precision factors
are done in high precision. Hence in exact arithmetic
\[
\vx_+ = \vx_c + \mm_{IR} \vx_c + {\hat \mU}^{-1} {\hat \ml}^{-1} \vb.
\]
So the residual update is
\[
\begin{array}{ll}
\vr_+ & = \vb - \ma \vx_c \\
\\
& = \vr_c - \ma {\hat \mU}^{-1} {\hat \ml}^{-1} \vr_c
\equiv \mm_{RES} \vr_c,
\end{array}
\]
where
\[
\mm_{RES} = \mi - \ma {\hat \mU}^{-1} {\hat \ml}^{-1}.
= {\hat \mU}^{-1} {\hat \ml}^{-1}
({\hat \ml} {\hat \mU} - \ma ).
\]
The analysis in the previous section implies that
\[
\| \mm_{RES} \|  \le \alpha,
\]
where 
\begeq
\label{eq:rate2}
\alpha = 
u_F ( \| \ma \| + 13 \sqrt{N} \kappa({\hat \ml})\kappa(\hat \mU) ).
\endeq
One could also use $\rho(\mm_{IR})$ for the convergence rate, but
we think \eqnok{rate2} is more illuminating. 

As is standard, when one computes a residual $\vr$ the computed
value ${\hat \vr}$ has an error \cite{higham}
\[
{\hat \vr} = \vr + \delta_\vr
\]
where
\[
\| \delta_\vr \| \le \gamma_N(u_W) (\| \ma \| \| \vx \| + \| \vb \|).
\]

We will do the analysis in terms of reduction in the residual norm. We
will then use that to estimate the limiting behavior of the error norm.
We will assume that $\alpha < 1$ and that the IR iteration is
bounded 
\[
\| \vx \| \le C \| \vx^* \|
\]
Hence
\begeq
\label{eq:residerr}
\| \delta_{\vr} \| \le \xi \equiv
\gamma_N(u_W) (C \| \ma \| \| \vx^* \| + \| \vb \|).
\endeq

We will analyze the progress of IR hile only considering the errors
in the the residual computation. So we compute
\[
\hat \vr_+ = \mm_{RES} {\hat \vr_c}  + \delta_{\vr_c}
\]
implying that
\[
\| \vr_+ \| \le \alpha \| \vr_c \| + (1 + \alpha ) \| \delta_{\vr_c} \| \le
\alpha \| \vr_c \| + (1 + \alpha ) \xi.
\]
Hence, for any $n \ge 0$
\[
\| \vr_{n+1} \| \le \alpha \| \vr_n \| + \frac{1 + \alpha}{1 - \alpha} \xi.
\]

So, the iteration will stagnate when 
\begeq
\label{eq:stagterm}
\| \vr \| \approx \frac{1 + \alpha}{1 - \alpha} \xi.
\endeq
When we terminate the iteration when $\| \vr \|/\| \vb \|$ is small
we are ignoring the $\| \ma \| \| \vx^* \|$ term in $\xi$, which is
one reason we must take watch for stagnation in our solver.

\subsubsection{Extended precision residuals}
\label{sec:extconv}

We consider the case $TR > TW$ and do interprecision transfers in
the triangular solves on the fly. In the language of 
\cite{CarsonHigham} the solver precision and working precision are the same.

The convergence result for this case (see Corollary 3.3, page A824, of
\cite{CarsonHigham} or also see \cite{moler}) is (stated very roughly) 

\begin{theorem}
\label{th:carsonh1}
Assume that $\| M_{IR} \|$ is sufficiently $< 1$. Then the error in IR
will be reduced until, at stagnation,
\begeq
\label{eq:irterm}
\frac{\| \vx - \vx^* \|_\infty}{\| \vx^* \|_\infty} 
\le 4 (N+1) \kappa_\infty(\ma) u_r + u_w.
\endeq
\end{theorem}
So, if $u_r \le u_w^2$ and $(N+1) \kappa_\infty(\ma) = O(u_w^{-1})$ then
\begeq
\label{eq:wilk1}
\frac{\| \vx - \vx^* \|_\infty}{\| \vx^* \|_\infty} = O(u_w),
\endeq
which is what we observed in the examples.
%
%





\appendix

\section{Interprecision Transfers: Part II}
\label{subsec:precision2}

In \cite{ctk:sirev20, ctk:fajulia, ctk:Newton3p}
we advocated LPS interprecision
with \eqnok{vdlow} rather than MPS. In this section we will look into
that more deeply and this discussion will reflect our recent experience
\cite{ctk:irnote}.
We will begin that investigation by comparing the cost
of triangular solves with the two approaches to interprecision transfer
to the cost of a single precision LU factorization. Since the triangular
solvers are $O(N^2)$ work and the factorization is $O(N^3)$ work,
the approach to interprecision transfer will matter less as the dimension
of the problem increases.

The test problem was $\ma \vx = \vb$ where the right side is $\ma$
applied to the vector with $1$ in each component. In this way we 
can compute error norms exactly.

\subsection{Double-Single IR}
\label{subsec:x6432}

In Table~\ref{tab:ip1} we report
timings from Julia's {\bf BenchmarkTools} package for double precision
matrix vector multiply (MV64),
single precision LU factorization (LU32) and three approaches
for using the factors to solve a linear system. HPS is the time for
a fully double precision triangular solved and MPS and LPS are the
mixed precision solve and the fully low precision solve using
\eqnok{vdlow} and \eqnok{promotevd}. IR will use a high precision
matrix vector multiply to compute the residual and a solve to 
compute the correction for each iteration. The low precision 
factorization is done only once.

\begin{table}[h]
\label{tab:ip1}
\caption{Timings for matrix-vector products and triangular 
solves vs factorizations: $\alpha=800$}
\centerline{
\begin{tabular}{lllllll} 
        N &     MV64 &     LU32 &      HPS &      MPS &      LPS & LU32/MPS \\ 
\hline 
  512  &  4.2e-05 &   1.2e-03  &  5.0e-05 &  1.0e-04  & 2.8e-05 & 1.2e+01   \\ 
 1024  &  8.2e-05 &   3.2e-03  &  1.9e-04 &  4.3e-04  & 1.0e-04 & 7.3e+00   \\ 
 2048  &  6.0e-04 &   1.4e-02  &  8.9e-04 &  2.9e-03  & 4.0e-04 & 4.8e+00   \\ 
 4096  &  1.9e-03 &   8.4e-02  &  4.8e-03 &  1.4e-02  & 2.2e-03 & 5.8e+00   \\ 
 8192  &  6.8e-03 &   5.8e-01  &  1.9e-02 &  5.8e-02  & 9.8e-03 & 1.0e+01   \\ 
\hline 
\end{tabular} 
}
\end{table}

The last column of the table is the ratio of timings for the low precision
factorization and the mixed precision solve. Keeping in mind that at least
two solves will be needed in IR, the table shows that MPS can be
a significant fraction of the cost of the solve for smaller problems and
that LPS is at least 4 times less costly. This is a compelling case
for using LPS in the case considered in this section, where high precision
is double and low precision is single, provided the performance of IR
is equally good.

If one is solving $\ma \vx = \vb$ for multiple right hand sides, as one
would do for nonlinear equations in many cases \cite{ctk:fajulia}, then
LPS is significantly faster for small and moderately large problems. For
example, for $N=4096$ the cost of MPS is roughly $15\%$ of the low precision
LU factorization, so if one does more than 6 solves with the same 
factorization, the solve cost would be more than the factorization cost. 
LPS is five times faster and we saw this effect while preparing 
\cite{ctk:fajulia} and we use that in our nonlinear solver package
\cite{ctk:siamfanl}.
The situation for IR is similar, but one must consider
the cost of the high precision matrix-vector multiply, which is about
the same as LPS. 

We make LPS the default for IR if high precision is double and low precision
is single. This decision is good for desktop computing. If low precision
is half, then the LPS vs MPS decision is different since the factorization
in half precision is so expensive.

Finally we mention a subtle programming
issue. We made Table~\ref{tab:ip1} with the
standard commands for matrix-vector multiply ($\ma * \vx$), factorization
{\bf lu}, and used $\backslash$ for the solve. Julia also offers non-allocating
versions of these functions. In Table~\ref{tab:ip1dot5} we show how using
those commands changes the results. We used {\bf mul!} for matrix-vector
multiply, {\bf lu!} for the factorization, and {\bf ldiv!} for the solve.

\begin{table}[h!]
\label{tab:ip1dot5}
\caption{Timings for non-allocating matrix-vector products and triangular 
solves vs factorizations: $\alpha=800$}
\centerline{
\begin{tabular}{lllllll} 
        N &     MV64 &     LU32 &      HPS &      MPS &      LPS & LU32/MPS \\ 
\hline 
  512  &  3.6e-05 &   9.1e-04  &  5.0e-05 &  4.8e-05  & 2.8e-05 & 1.9e+01   \\ 
 1024  &  9.0e-05 &   2.7e-03  &  1.9e-04 &  1.8e-04  & 1.0e-04 & 1.5e+01   \\ 
 2048  &  6.2e-04 &   1.3e-02  &  8.9e-04 &  7.3e-04  & 3.9e-04 & 1.8e+01   \\ 
 4096  &  2.2e-03 &   8.0e-02  &  4.8e-03 &  3.3e-03  & 2.3e-03 & 2.4e+01   \\ 
 8192  &  6.5e-03 &   5.7e-01  &  2.1e-02 &  1.5e-02  & 1.0e-02 & 3.9e+01   \\ 
\hline 
\end{tabular} 
}
\end{table}

So, while LPS still may make sense for small problems if high precision is
double and low precision is single, 
the case for using it is weaker if one uses non-allocating 
matrix-vector multiplies and solves. We do that in 
{\bf MultiPrecisionArrays.jl}.

\subsection{Accuracy of MPS vs LPS}
\label{sec:mpsacc}

Since MPS does the triangular solves in high precision, one should expect
that the results will be more accurate and that the improved accuracy
might enable the IR loop to terminate earlier \cite{CarsonHigham}.
We should be able to see
that by timing the IR loop after computing the factorization. One should
also verify that the residual norms are equally good.

We will conclude this section with two final tables for the results of IR.
We compare the well
conditioned case ($\alpha=1$) and the ill-conditioned case ($\alpha=800$)
for a few values of $N$. We will look at residual and error norms
for both approaches to interprecision transfer. The conclusion is that
if high precision is double and low is single, the two approaches give
equally good results.  

The columns of the tables are the dimensions, the 
$\ell^\infty$ relative error norms for both
LP and MP interprecision transfers (ELP and EMP) and the corresponding
relative residual norms (RLP and RMP). 

The results for $\alpha=1$ took 5 IR iterations for all cases. As expected
the LPS iteration was faster than MPS. 
However,
for the ill-conditioned $\alpha=800$ case, MPS took one fewer iteration 
(5 vs 6)
than EPS
for all but the smallest problem. 
Even so, the overall solve
times were essentially the same. 

\begin{table}[h]
\label{tab:ip2}
\caption{Error and Residual norms: $\alpha=1$}
\begin{center}
\begin{tabular}{lllllll} 
        N &      ELP &      EMP &      RLP &      RMP &      TLP &      TMP \\ 
\hline 
  512  &  4.4e-16 &   5.6e-16  &  3.9e-16 &  3.9e-16  & 3.1e-04 & 3.9e-04   \\ 
 1024  &  6.7e-16 &   4.4e-16  &  3.9e-16 &  3.9e-16  & 1.1e-03 & 1.5e-03   \\ 
 2048  &  5.6e-16 &   4.4e-16  &  3.9e-16 &  3.9e-16  & 5.4e-03 & 6.2e-03   \\ 
 4096  &  1.1e-15 &   1.1e-15  &  7.9e-16 &  7.9e-16  & 1.9e-02 & 2.5e-02   \\ 
 8192  &  8.9e-16 &   6.7e-16  &  7.9e-16 &  5.9e-16  & 6.9e-02 & 9.3e-02   \\ 
\hline 
\end{tabular} 
\end{center}
\end{table}

\begin{table}[h]
\label{tab:ip3}
\caption{Error and Residual norms: $\alpha=800$}
\begin{center}
\begin{tabular}{lllllll} 
N &      ELP &      EMP &      RLP &      RMP &      TLP &      TMP \\ 
\hline 
  512  &  6.3e-13 &   6.2e-13  &  2.1e-15 &  1.8e-15  & 3.0e-04 & 3.8e-04   \\ 
 1024  &  9.6e-13 &   1.1e-12  &  3.4e-15 &  4.8e-15  & 1.4e-03 & 1.5e-03   \\ 
 2048  &  1.0e-12 &   1.2e-12  &  5.1e-15 &  4.5e-15  & 6.5e-03 & 7.1e-03   \\ 
 4096  &  2.1e-12 &   2.1e-12  &  6.6e-15 &  7.5e-15  & 2.6e-02 & 2.4e-02   \\ 
 8192  &  3.3e-12 &   3.2e-12  &  9.0e-15 &  1.0e-14  & 9.1e-02 & 8.7e-02   \\ 
\hline 
\end{tabular} 
\end{center}
\end{table}



\bibliographystyle{siamplain}
\bibliography{MPArray}


\end{document}